\documentstyle[12pt]{article}

\def\be{\begin{equation}}
\def\ee{\end{equation}}
\def\bea{\begin{eqnarray}}
\def\eea{\end{eqnarray}}
\def\ba{\begin{array}}

\def\ea{\end{array}}
\begin{document}
\vskip .2in
\hfill
\vbox{
    \halign{#\hfil         \cr
           IPM-95-100      \cr
           August 1995   \cr
           } 
      }  
\vskip 3cm
\centerline{\bf \Large  Generalization of the $h$-Deformation to Higher
Dimensions\footnote{To be appear in J. Phys. A}}
\centerline{M. Alishahiha}
\centerline{\it Institute for Studies in Theoretical Physics and Mathematics,}
\centerline{\it  P.O.Box 19395-5746, Tehran, Iran }
\centerline{\it  Department of Physics, Sharif University of Technology,}
\centerline{\it  P.O.Box 11365-9161, Tehran, Iran }
\centerline{\it e-mail: alishah@netware2.ipm.ac.ir}
\vskip 3cm
\begin{abstract}
In this article we construct $GL_{h}(3)$ from $GL_{q}(3)$ by a singular map.
We show that there exist two singular maps which map $GL_{q}(3)$ to new
quantum groups. We also construct their $R$-matrices and will show
although the maps are singular but their $R$-matrices are not. Then we
generalize these singular maps to the case $GL(N)$ and for $C_{n}$ series.
\end{abstract}
\newpage

There exist two types of $ SL(2) $ quantum groups. One is the standard
$ SL_{q}(2)$, another one is the Jordanian quantum group which is also
called the $h$-deformation of $SL(2) $. Quantum matrices in two dimensions,
admitting left and right quantum spaces, are classified [1].
One is the $q$-deformation of $ GL(2) $, the other is the $h$-deformation.
 The $q$-deformation of $ GL(N) $ has been studied extensively but in the
 literature only the two dimensional case of $h$-deformation has been
studied.[2-7]

In ref. 8 it is shown that $ GL_{h}(2) $ can be obtained from $ GL_{q}(2) $
by a singular limit of a similarity transformation. We will show this method
can be used successfully, for construction of $GL_{h}(N)$. In other words,
at first we will consider the $GL(3)$, and introduce two singular maps
which convert $GL_{q}(3)$ to $GL_{h}(3)$. Then we generalize one of the
singular maps to $N$-dimensional case. We will use $R$-matrix of $GL_{q}(N)$
which by this map, results to a new $R$-matrix. Also, by this map one can
obtain $h$-deformation of $C_{n}$ series, but can not for $B_{n}$ and $D_{n}$
 series.

In this article we denote $q$-deformed objects by primed quantities. Unprimed
quantities represent transformed objects.

Consider Manin's $q$-plane with the following quadratic relation between
coordinates.
\be
x'y'=q y'x'.
\ee
By the following linear transformation:
\bea
\pmatrix{x' \cr y'} &=&\pmatrix{1 & {{h}\over{q-1}} \cr 0 & 1} \pmatrix{x \cr
y}
\eea
the relation $(1)$ changes to $ x y - q y x = h y^{2} $. For the case of
$ q=1 $, one get the relation of two dimensional $h$-plane. In fact $g$ itself
is
 singular in the $ q=1 $ case, but the resulting relation for the plane is
 non-singular.

The above linear transformation on the plane induces the following similarity
transformation on the $R$-matrix of $ GL_{q}(2) $.
\be
R_{h}=\lim_{q\rightarrow 1 }{(g \otimes g)}^{-1} R_{q} (g \otimes g).
\ee

Although the above map is singular, the resulting $R$-matrix is non-singular
and is the well known $R$-matrix of $ GL_{h}(2) $.

Now consider 3-dimensional Manin's quantum space:
\be\ba{ll}
 x'_{i} x'_{j}=q x'_{j} x'_{i} &  i<j,
 \ea\ee
 and consider the following linear transformation:
\be
X=g^{-1}X',
\ee
where
\be
g=\pmatrix{\lambda_{1} & \alpha & \beta \cr 0 & \lambda_{2} & \gamma \cr
0 & 0 & \lambda_{3}}.
\ee

Here $\alpha, \beta$ and $\gamma$ are parameters which can be singular at
$q=1$. So they can be written as ${1\over{f(q)}}$ where $f(1)=0$.
The Taylor expansion of $f(q)$ about $q=1$ is $f(q)={1 \over{h}}(q-1)+
O((q-1)^{2})$. We need only the first term, because we are only interested
in the behaviour of $f(q)$ in the neighbourhood of $q=1$. The
coefficient of
first term in the Taylor expansion, $h$, plays the role of the deformation
parameter for the new quantum group. The $\lambda_{i}$s can be made equal to 1
by rescaling.

To obtain $\alpha, \beta$ and $\gamma$ we should apply this map to the
$q$-deformed plane and its dual, and require that the mapped plane and its
dual be non-singular at $q=1$. The following are the only singular maps
satisfying this condition:

\be\ba{lll}
g_{1}=\pmatrix{1 &{h \over{q-1}} & \beta \cr 0 & 1 & 0 \cr 0 & 0 & 1}, &
g_{2}=\pmatrix{1 & \alpha & \beta \cr 0 & 1 & {h \over{q-1}} \cr 0 & 0 & 1}, &
g_{3}=\pmatrix{1 & \alpha &{h \over{q-1}} \cr 0 & 1 & \gamma \cr 0 & 0 & 1}.
\ea\ee

Here $\alpha,\ \beta$ and $\gamma$ (in $g_{1},\ g_{2},\ g_{3},$\ ) are
non-singular
parameters. Note that the $R$-matrices obtained from these maps,
solve the quantum Yang-Baxter equation and are non-singular for $q=1$.

Let us denote the dependence of $g_{1}, g_{2}$ and $g_{3}$ on parameters
explicitly:
\be\ba{lll}
g_{1} := g_{1}({h\over(q-1)},\beta), &
g_{2} := g_{2}({h\over(q-1)},\alpha,\beta),&
g_{3} := g_{3}({h\over(q-1)},\alpha,\gamma).
\ea\ee
It is easy to show that:
\bea
g_{1}({h\over(q-1)},\beta) g_{1}(0,-\beta)&=& g_{1}({h\over(q-1)},0) \cr
g_{2}({h\over(q-1)},\alpha,\beta) g_{2}(0,-\alpha,-\beta)&=&
g_{2}({h\over(q-1)},0,0) \cr
g_{3}({h\over(q-1)},\alpha,\gamma)g_{3}(\alpha \gamma,-\alpha,-\gamma) &=&
g_{3}({h\over(q-1)},0,0).
\eea
so all non-singular parameters in the above matrices can be set to zero.
Moreover the
$R$-matrices $ R(g_{1})$ and $R(g_{2})$ which are obtained by formula (3)
using $g_{1}({h\over(q-1)},0)$ and $g_{2}({h\over(q-1)},0,0)$ respectively,
are equivalent, because:
\be
{(s \otimes s)}^{-1} R(g_{2}) (s \otimes s)=R(g_{1}).
\ee
where $s=e_{13}+e_{21}+e_{32}$. So, there are only two independent cases.
The $R$-matrices corresponding to these transformations are non-singular
and have been first obtained by Hietarinta [9]. The first case ( the trivial
case) is $\beta=0$ in $g_{1}$ (or $\alpha=\beta=0$ in $g_{2}$) and the second
case is $\alpha=\gamma=0$ in $g_{3}$. The $h$-deformed quantum plane and its
dual and $R$-matrices corrsponding to these cases are:

{\bf{First case}}
\be\ba{ll}
[x_{1},x_{2}]=hx_{2}^{2}, & \eta_{3}^{2}=\eta_{2}^{2}=\{\eta_{1},\eta_{2}\}=0,
\cr
[x_{1},x_{3}]=0, & \{\eta_{2},\eta_{3}\}=\{\eta_{1},\eta_{3}\}=0, \cr
[x_{2},x_{3}]=0, & \eta_{1}^{2}=-h\eta_{2}\eta_{1}.
\ea\ee
and the non-zero elements of $R$-matrix except for $R_{ijij}=1$ are:
\bea
R_{1121}=R_{2122}&=&-R_{1112}=-R_{1222}=h, \cr
R_{1122}&=&h^{2}.
\eea

{\bf{Second case}}
\be\ba{ll}
[x_{1},x_{2}]=2hx_{3}x_{2}, & \{\eta_{1},\eta_{2}\}=-2h\eta_{3}\eta_{2},
\cr [x_{1},x_{3}]=hx_{3}^{2}, & \eta_{1}^{2}=-h\eta_{3}\eta_{1}, \cr
[x_{2},x_{3}]=0, & \eta_{3}^{2}=\eta_{2}^{2}=\{\eta_{1},\eta_{3}\}=\{\eta_{2},
\eta_{3}\}=0.
\ea\ee
and the non-zero elements of $R$-matrix except for $R_{ijij}=1$ are:
\be\ba{ll}
R_{1113} = R_{1333}=-h, & R_{1131} = R_{3133}=h, \cr
R_{2132} = -R_{1223}=2h & R_{1133} $ = $h^{2}.
\ea\ee

A linear transformation on the plane induces a similarity transformation on the
quantum matrices acting upon it.
\be
 M'= g M g^{-1},
\ee

The algebra of functions, $ GL_{q}(3) $, is obtained from the following
relations:
\be
R' M'_{1} M'_{2} = M'_{2} M'_{1} R'.
\ee

Applying transformation (15) one easily obtains for the case of $ q=1
$. \be
R M_{1} M_{2}= M_{2} M_{1} R.
\ee

So the entries of the transformed quantum matrix $ M $ fulfill the commutation
relations of the $ GL_{h}(3)$, for both $g$'s. It is easy to show that the
$h$-deformed determinant is central, so it can be set to 1. A quantum group's
differential structure is completely determined by $R$-matrix [10].
One therefore expects that by these similarity transformations the differential
structure of the $h$-deformaton be obtained from that of the $q$-deformation.
\bea
M_{2}dM_{1} - R_{12}dM_{1}M_{2}R_{21}&=&0, \cr
dM_{2}dM_{1} + R_{12}dM_{1}dM_{2}R_{21}&=&0.
\eea

Now, it is obvious that, defining $dM:=g^{-1}dM g$ and using the above
relations
the differential of $GL_{h}(3)$ can be easily obtained from the corresponding
differential structure of $GL_{q}(3)$.

For the higher dimensions, there are several generalizations which depend on
the position of singularity in $g$. For example we consider the following
generalization: \be
g=\sum_{i=1}^{N} e_{ii} + {{h} \over{q-1}} e_{1N}
\ee

The general aspect of the contraction for arbitrary $N$ can be obtained from
this
simple map. By inserting this map in (3)  we will obtain the general form of
the $h$-deformed $R$-matrix, which solves the quantum Yang-Baxter equation.

{\bf{1- The series $ A_{n-1} $ }}

After applying this singular map, the corresponding $h$-deformed $R$-matrix
will become:
\bea
R_{h}&=&\sum_{i,j=1}^{N} e_{ii} \otimes e_{jj}
+ 2h \sum_{i>1}^{N-1} (e_{1i} \otimes e_{iN} - e_{iN} \otimes e_{1i}) \cr
&-&h(e_{1N} \otimes e_{NN} - e_{NN} \otimes e_{1N}) -h(e_{11} \otimes e_{1N} -
e_{1N} \otimes e_{11})  \cr&+&h^{2}(e_{1N} \otimes e_{1N}).
\eea

Consider $N$-dimensional $q$-deformed quantum space
\be\ba{ll}
x'_{i} x'_{j} = q x'_{j} x'_{i} & i<j.
\ea\ee

Assume the following linear singular transformation
\be
x'_{i} = g_{ij} x_{j}.
\ee

By the above transformation and in the $ q=1$ case we obtain the $h$-deformed
quantum plane as follows:
\be\ba{ll}
x_{i} x_{j} = x_{j} x_{i}  &  1<i<j\leq N, \cr
[x_{1},x_{j}]=2h x_{N} x_{j}, & [x_{1},x_{N}]= h (x_{N})^{2}. \cr
\ea\ee

{\bf{2- The series $ B_{n}, C_{n} $ and $ D_{n}$ }}

The corresponding $q$-deformed $R$-matrix has order $ N^2 \times N^2 $, where $
N=2n+1$
 for $ B_{n} $ and $ N=2n $ for $ D_{n} $ and $ C_{n} $ and it is given by
[11]:
\bea
R_{q}&=&q \sum_{i \neq i'}^{N}e_{ii}\otimes e_{ii} +e_{{{N+1}\over{2}}
{{N+1}\over{2}}}
\otimes e_{{{N+1}\over{2}} {{N+1}\over{2}}}+\sum_{i \neq j,j'}^{N}
e_{ii}\otimes e_{jj}\cr
&+&q^{-1} \sum_{i \neq i'}^{N} e_{i'i'} \otimes e_{ii} +(q-q^{-1}) \sum_{i>j}^
{N} e_{ij} \otimes e_{ji} \cr
&-&(q-q^{-1}) \sum_{i>j}^{N} q^{\rho_{i}-\rho_{j}} \epsilon_{i} \epsilon_{j}
e_{ij} \otimes e_{i'j'}.
\eea

The second term is present only for the series $ B_{n} $. Here $ i'=N+1-i,
j'=N+
1-j, \epsilon_{i}=1, i=1,...,N $ for the series $ B_{n} $ and  $ D_{n}$,
$ \epsilon_{i}=1, i=1,...,{{N}\over{2}}$, $\epsilon_{i}=-1,
i={{N}\over{2}}+1,...,N $
for the series $ C_{n} $ and $ (\rho_{1},...,\rho_{N}) $ is:
\be\ba{ll}
({n-{1}\over{2}},...,{{1}\over{2}},0,{{1}\over{2}},...,-n+{{1}\over{2}}) & for
B_{n} \cr
(n,n-1,...,1,-1,...,-n) & for C_{n} \cr
(n-1,...,1,0,0,-1,...,-n+1) & for D_{n}
\ea\ee

By inserting this $R$-matrix in (3), the coefficient of $ e_{1N} \otimes e_{1N}
$
will become:
\be
{{h^{2}}\over{q-1}} (q^{-1}+1)(1+\epsilon_{N} q^{\rho_{N}-\rho_{1}}).
\ee

This expression is non-singular only when $\epsilon_{N}=-1$ and for $ q=1 $
it is equal to $ 2 N h^{2} $. We thus see that only the $ C_{n} $ series
remains non-singular. The corresponding $h$-deformed $R$-matrix is:
\bea
R_{h}&=&\sum_{i,j=1}^{N} e_{ii} \otimes e_{jj} + 2N h^{2} e_{1N}\otimes
e_{1N}\cr
&-&2h \sum_{i=2}^{N} e_{1i} \otimes e_{iN}+\epsilon_{i} e_{iN}\otimes e_{i'N}
\cr
&+&2h \sum_{i=1}^{N-1}e_{iN}\otimes e_{1i} -\epsilon_{i} e_{1i}\otimes e_{1i'}.
\eea

So by this method we can obtain $ SP_{h}(2n) $. The algebra $ SP^{2n}_{q}(c) $
with generators $ x'_{1},...,x'_{2n} $ and relations
\be
R'_{q}(x'\otimes x')= q x'\otimes x',
\ee
is called the algebra of functions on quantum $2n$-dimensional symplectic
space.
After applying the singular transformation (19) to (28) we obtain the relations
between the generators of $ SP_{h}^{2n}(c) $:
\be\ba{lll}
x_{i}x_{j}=x_{j}x_{i},\   &  1<i<j\leq N,\  & j \neq j',
\ea\ee
\be\ba{ll}
x_{1}x_{j}=x_{j}x_{1}+2h x_{N}x_{j},\  &  j \neq N, \cr
x_{i'}x_{i}=x_{i}x_{i'}+2h \epsilon_{i'} x^{2}_{N},\  & 1<i<i'\leq N.
\ea\ee

In $ SP^{2n}_{q}(c)$ the equality $ x'^{t}C'x'=0 $ holds. By applying the
singular map (29), $ C'$ transforms to $ C=g^{t} C' g $, where $C$ is given by:
\be
C=\sum_{i=1}^{N} \epsilon_{i} e_{ii'} -Nh e_{NN}.
\ee

The Quantum group $ SP_{q}(2n) $ acts on $ SP^{2n}_{q}(c) $ and preserves
$ x'^{t} C' x'=0 $, so we have:
\be
M'^{t} C' M'=C',
\ee
on the other hand:
\bea
M=gM'g^{-1},\        & M^{t}=(g^{-1})^{t} M'^{t} g^{t}.
\eea
It follows that:
\be
M^{t} C M= C,
\ee

So we conclude that the quantum group $SP_{h}(2n)$ acts on $SP_{h}^{2n}(c)$
and preserves $ x^{t} C x=0 $.
It is interesting to note that the expression $ x'^{t} C' x'$, which should
 be equal to 1 for $ SO(2n) $ and $ SO(2n+1) $ ($ B_{n} $ and $ D_{n} $
series),
is singular. So we cannot obtain the $h$-deformation of $ B_{n} $ and $ D_{n} $
series by contraction of the $q$-deformation, at least by this form (upper
triangular matrix) of singular transformation ($g$).

One of the interesting problems is to construct $ U_{h}(gl(3)) $, and its
generaliztion to higher dimensions.

{\bf{Acknowledgments}}

I would like to thank A. Aghamohammadi for drawing my attention to this
problem,
and V. Karimipour, A. Shariati, and M. Khorami for valuable discussions. I
would also like to thank referee for his(her) useful comments, specially for
the discussion on classifying the non-equivalant singular transformation.

{\bf{References}}

1. H. Ewen, O. Ogievetsky and J. Wess, Lett. Math. Phy. 22,297 (1991)

2. E. Domidov, Yu. I. Manin, E. E. Mukhin and D. V. Zhdanovich, Prog. Theor.

Phys. Suppl. 102, 203 (1990)

3. S. Zakrzewki, Lett. Math. Phys. 22, 287 (1991)

4. Ch. Ohn, Lett. Math. Phy. 25, 891 (1992)

5. V. Karimipour, Lett. Math. Phys. 30, 87 (1994)

6. B. A. Kupershmidt, J. Phys. A 25, L123 (1992)

7. A. Aghamohammadi, Mod. Phy. Lett. A 8 (1993)

8. A. Aghamohammadi, M. Khorrami and A. Shariati, J. Phys. A. 28, L225 (1995)

9. J. Hietarinta, J. Phy. A. 26, 7077 (1993)

10. A. Sudbery, Phys. Lett. B 284, 61 (1992)

11. N. Yu. Reshetikhin, L. A. Takhtadzhyan, L. D. Faddeev, Len. Math. J. vol. 1

No. 1 (1990)

\end{document}